\documentclass[11pt]{amsart}
\usepackage{amsmath}
\usepackage{xypic}
\usepackage{amssymb}
\usepackage{hyperref}

\newtheorem{theorem}{Theorem}[section]
\newtheorem{proposition}[theorem]{Proposition}
\newtheorem{lemma}[theorem]{Lemma}
\newtheorem{definition}[theorem]{Definition}
\newtheorem{examples}[theorem]{Examples}

\newtheorem{corollary}[theorem]{Corollary}
\newtheorem{remark}[theorem]{Remark}
\newtheorem{example}[theorem]{Example}

\begin{document}

\def\be{\begin{equation}}
\def\ee{\end{equation}}

\def\ra#1{\mathop{\longrightarrow}\limits^{#1}}
\def\rab#1{\mathop{\longrightarrow}\limits_{#1}}
\def\lab#1{\mathop{\longleftarrow}\limits_{#1}}
\def\la#1{\mathop{\longleftarrow}\limits^{#1}}
\def\da#1{\downarrow\rlap{$\vcenter{\hbox{$\scriptstyle#1$}}$}}
\def\ua#1{\uparrow\rlap{$\vcenter{\hbox{$\scriptstyle#1$}}$}}
\def\sea#1{\mathop{\searrow}\limits^{#1}}
\def\nea#1{\mathop{\nearrow}\limits^{#1}}
\def\swa#1{\mathop{\swarrow}\limits^{#1}}
\def\dirlim#1{\lim_{\textstyle{\mathop{\rightarrow}\limits_{#1}}}}
\def\invlim#1{\lim_{\textstyle{\mathop{\leftarrow}\limits_{#1}}}}

\def\AA{\mathbb{A}}
\def\DD{\mathbb{D}}
\def\QQ{\mathbb{Q}}
\def\ZZ{\mathbb{Z}}
\def\NN{\mathbb{N}}
\def\GG{\mathbb{G}}
\def\WW{\mathbb{W}}
\def\RR{\mathbb{R}}
\def\CC{\mathbb{C}}
\def\FF{\mathbb{F}}
\def\HH{\mathbb{H}}
\def\LL{\mathbb{L}}

\def\qed{$\,\square$}

\def\tensor{\otimes}
\def\dsum{\oplus}
\def\Tensor{\mathop\bigotimes}
\def\Dsum{\mathop\bigoplus}
\def\img{\mathrm{im}}

\def\bR{\mathbf{R}}

\def\cW{\mathcal{W}}
\def\cC{\mathcal{C}}
\def\cF{\mathcal{F}}
\def\cG{\mathcal{G}}
\def\cO{\mathcal{O}}
\def\cL{\mathcal{L}}
\def\cM{\mathcal{M}}
\def\cN{\mathcal{N}}

\def\bSet{\mathbf{Set}}
\def\bMod{\mathbf{Mod}}
\def\bAb{\mathbf{Ab}}
\def\bC{\mathbf{C}}
\def\bA{\mathbf{A}}
\def\bMag{\mathbf{Mag}}
\def\bV{\mathbf{V}}
\def\bAss{\mathbf{Ass}}
\def\bCom{\mathbf{Com}}
\def\bAlt{\mathbf{Alt}}
\def\bAAlt{\mathbf{AAlt}}
\def\bDelta{\mathbf{\Delta}}
\def\bSigma{\mathbf{\Sigma}}
\def\bOrd{\mathbf{Ord}}
\def\bLie{\mathbf{Lie}}
\def\bNov{\mathbf{Nov}}
\def\bVin{\mathbf{Vin}}
\def\bAlg{\mathbf{Alg}}
\def\bCA{\mathbf{CAlg}}
\def\bLeib{\mathbf{Leib}}
\def\bLev{\mathbf{Lev}}
\def\bRMod{\mathbf{RMod}}
\def\bLMod{\mathbf{LMod}}

\def\Hom{\mathrm{Hom}}
\def\Map{\mathrm{Map}}
\def\c{\mathrm{c}}
\def\Ab{\mathrm{Ab}}
\def\End{\mathrm{End}}
\def\pr{\mathrm{pr}}
\def\Tens{\mathrm{Tens}}
\def\gr{\mathrm{gr}}
\def\Ext{\mathrm{Ext}}
\def\sgn{\mathrm{sgn}}
\def\res{\mathrm{res}}
\def\rank{\mathrm{rank}}

\def\intersect{\mathop\bigcap}
\def\iso{\cong}
\def\cotensor{\square}

\def\ext{{e}}

\author{Nishant Dhankhar}
\address{Department of Mathematics, Massachusetts Institute of Technology,
Cambridge, MA 02139}
\email{dhankhar@mit.edu}

\author{Haynes Miller}
\address{Department of Mathematics, Massachusetts Institute of Technology,
Cambridge, MA 02139}
\email{hrm@math.mit.edu}

\author{Ali Tahboub}
\address{Department of Mathematics, Birzeit University,
Ramallah, Palestine}
\email{aliraeedtahbub@hotmail.com}

\author{Victor Yin}
\address{Department of Mathematics, Massachusetts Institute of Technology,
Cambridge, MA 02139}
\email{victor05@mit.edu}

\title
{Beck modules and alternative algebras}

\subjclass[2010]{17A30,17.04}

\begin{abstract}
We set out the general theory of ``Beck modules'' in a variety of
binary algebras and
describe them as modules over suitable ``universal enveloping'' unital
associative algebras. We develop a theory of ``noncommutative partial 
differentiation'' to pass from the equations of the variety to relations in 
a universal enveloping algebra. We pay particular attention to the
case of alternative algebras, defined by 
a restricted associative law, and determine the Poincar\'e polynomial
of the universal enveloping algebra in the abelian case. 
\end{abstract}

\setcounter{equation}{0}

\maketitle

\section{Introduction}
The notion of a ``module'' occupies an important place in the study
of general algebraic systems. Most of these diverse notions are united
under the theory of ``Beck modules.'' Given an object $A$ in any category
$\bC$, one may consider the ``slice category'' $\bC/A$
of objects in $\bC$ equipped with a map to $A$.
A Beck module for $A$ is then an abelian group object in $\bC/A$.
If $\bC$ is the category of commutative rings, for example,
a Beck module for $A$ is simply an $A$-module, while if
$\bC $ is the category of associative algebras, a Beck module for $A$ is an
$A$-bimodule. Many other examples occur in the literature: Leibniz
algebras \cite{loday-pirashvili}, $\lambda$-rings
\cite{hesselholt}, divided power rings \cite{dokas}, \ldots.

This definition occurs in the thesis \cite{beck}
of Jonathan Beck written under the direction
of Samuel Eilenberg. Eilenberg himself had discussed 
such objects in \cite{eilenberg}, at least in the linear
context, as the kernel of a ``square zero extension.''
These kernels were understood to constitute ``representations'' of the 
algebra, and this structure was made explicit in various cases.

We review below the context of a ``variety'' $\bV$ of algebras over a
commutative ring $K$. In this case, for every $\bV$-algebra $A$ the category
$\bMod_A$ of Beck $A$-modules is an abelian category with a single
projective generator. As a result, the category $\bMod_A$ is equivalent
to the category of right modules over a canonical unital
associative $K$-algebra $U_\bV(A)$, the ``universal enveloping
algebra'' for $A$. 

This raises the question of identifying the structure of $U_\bV(A)$ for
various varieties $\bV$ and $\bV$-algebras $A$. Left and right multiplication
determine a $K$-module map $A\oplus A\to U_\bV(A)$, and hence a surjection of
associative unital $K$-algebras $\Tens_K(A\oplus A)\to U_\bV(A)$
(cf. \cite{loday-pirashvili}). 
Each defining equation determines a generator of the kernel of this map,
by a process of ``noncommutative differentiation'' that we describe in detail. 

We review some of the standard examples, and
then focus on a somewhat less standard one, the variety of ``alternative
algebras'' over $K$. This example has been considered before,
but even over a field basic features of the universal
enveloping algebra for an alternative algebra, such as its dimension,
have remained obscure. In 1954, Nathan Jacobson \cite{jacobson} wrote
``The introduction of the universal associative algebras for the
birepresentations [his term for Beck modules] enables one to split
the representation problem into
two parts: (1) determination of the structure of $U(A)$, (2) representation
theory for the associative algebra $U(A)$. In practice, however, it seems
to be difficult to treat (1) as a separate problem. Only in some
special cases is it feasible to attack this directly.''
Richard Schafer \cite{schafer} observed in 1966 that if $\bV$ is the
variety of alternative $K$-algebras, with $K$ a field, and $\dim_KA=n$,
then $\dim_KU(A)\leq4^n$. But the precise dimension, even the case
of ``abelian'' alternative algebras -- those with trivial product --
over a field, has eluded analysis.

In that case, the universal enveloping
algebra admits a natural grading, by word length, or, as we call it,
by weight. Write $K^n$ for the free $K$-module on $n$ generators, regarded as
an abelian alternative $K$-algebra. A signal 
result of this paper is the description of an explicit basis
for the universal enveloping algebra of $K^n$ when $K$ is a field.

\begin{theorem} Let $K$ be any commutative ring. For each $n$ and $k$,
the $K$-module $U(K^n)_k$ is free, and
\[
\rank_K U(K^n)_k=
\begin{cases} 1 & \text{if}\,\,k=0\\
2n & \text{if}\,\, k=1\\
\frac{3n^2-n}{2} & \text{if}\,\, k=2\\
2\textstyle{\binom{n}{k}} & \text{if}\,\,k\geq3\,.
\end{cases}
\]
\end{theorem}
In particular,
\[
U(K^n)_k=0\quad\text{for}\quad k>n\,.
\]

These calculations show that at least in the abelian case, the growth
rate of the universal enveloping algebra is indeed exponential in the
dimension of the algebra, but much
slower than the upper bound observed by Schafer.

Our tool is the theory of Gr\"obner bases for noncommutative graded algebras.
We employ hand calculation and Python to determine a 
Gr\"obner basis for the ideal of relations defining $U(K^n)$, with $K$ 
any prime field, for $n\leq 5$. The structure of this basis for
these small values of $n$ turns out to imply that a basis with the same
structure exists for all $n$. The set of normal monomials with respect
to this basis (which projects to a basis for $U(K^n)$) 
is then easy to determine. A base-change result then shows that 
$U(K^n)$ is a free $K$-module of the same rank for any commutative base
ring $K$. 

Returning for a moment to the situation of a general variety $\bV$, 
any right module
$M$ over the universal enveloping algebra $U_\bV(A)$ of a $\bV$-algebra $A$ 
determines a new $\bV$-algebra, the abelian object over $A$ corresponding 
to $M$. It has the form $E=A\oplus M$ as a $K$-module. If $A$ has trivial 
multiplication, this $K$-algebra has the property that $(uv)(xy)=0$ for
all $u,v,x,y\in E$. In the absence of associativity, this ``solvability''
condition does not imply a nilpotence condition, which would say that every
sufficiently long product, arbitrarily bracketed, vanishes. Taking $M$ to be 
the free right $U_\bV(A)$-module on one generator provides an especially
interesting case. With $\bV=\bAlt_K$ and $A=K^n$ with trivial product,
our work provides a nontrivial product
of length $n$ in $E$. We can take $n=\infty$, and obtain an example of an
alternative algebra of ``solvability index'' 2 that is not nilpotent.
A referee has pointed out to us that this is exactly the example published 
in Russian in 1960 by G. V. Dorofeev
\cite{dorofeev}, described also in the book \cite[\S6.2]{zhevlakov-et-al}.
We feel that setting this example in the broader context of Beck modules 
and universal enveloping algebras, along with the formalism of 
noncommutative partial differentiation, has independent value, as does the 
structured computation using Gr{\"o}bner bases. 

After a review of the theory of varieties of algebras in 
\S\ref{sec-varieties}, we describe in \S\ref{sec-beck} 
the theory of Beck modules in this generality, and the
corresponding universal enveloping algebras. In \S\ref{sec-partial-diff} 
we discuss the form of non-commutative differentiation that leads from
the equations in a variety $\bV$ to the relations in a universal enveloping 
algebra for a $\bV$-algebra. In \S\ref{sec-aa-univ} 
we specialize to the case of trivial alternative algebras, and in the
following section we present a few other results about universal enveloping
algebras of alternative algebras. Finally,
in an Appendix, we review some of the essential features of the theory of 
Gr\"obner bases.

This work is intended as a first step in the study of the Quillen homology
and cohomology of algebraic systems such as alternative algebras.

\medskip
\noindent
{\bf Acknowledgements.}
This work was carried out under the auspices of a program, supported by
MIT's Jameel World Education Laboratory, designed to foster
collaborative research projects involving students from MIT and Palestinian
universities.
We acknowledge with thanks the contributions made by early participants
in this program -- Mohammad Damaj and Nasief Khlaif of Birzeit University, 
Hadeel AbuTabeekh of An-Najah National University, and Heidi Lei of MIT --
as well as the 
support of Palestinian faculty -- Reema Sbeih and Mohammad Saleh at Birzeit
and Khalid Adarbeh and Muath Karaki at NNU. 
The first author acknowledges support by the MIT UROP office. 
We thank a referees for their interest and for pointing out several relevant 
references.

\section{Varieties of algebras}
\label{sec-varieties}

We will work with algebras defined by a product operation,
though much of this work can be carried out
in much greater generality. Following the lead of 
Bourbaki \cite[\S7.1]{bourbaki-algebra-1}, we make the following definition.

\begin{definition} A {\em magma} is a set $X$ with a binary operation
$X\times X\to X$ (written as juxtaposition). 
\end{definition}

Write $\bMag$ for the category of magmas. 

We will also restrict our attention to linear examples,
and work over a commutative ring $K$.
So a {\em magmatic $K$-algebra} (or just $K$-{\em algebra}) 
is a $K$-module $A$ equipped with a $K$-bilinear product 
inducing a linear map $A\tensor A\to A$, written as juxtaposition. 

Magmatic $K$-algebras constitute the objects in a category $\bMag_K$.
The forgetful functor to sets has a left adjoint $Mag$.
The free magma $Mag(S)$ generated
by a set $S$ is the set of bracketed strings of elements of $S$; see 
\cite[\S7.1]{bourbaki-algebra-1}. The free magmatic $K$-algebra on a set $S$ is
the free $K$-module generated by $Mag(S)$: $Mag_K(S)=KMag(S)$. 

We can adjoin axioms using the following device \cite{adamek}. 
An {\em equation} is an element $\omega$ of the free magmatic $K$-algebra on a 
finite set $S$, which we may denote by $S(\omega)$ if there are several
equations in play. Given a 
magmatic $K$-algebra $A$, we will say that an equation $\omega\in Mag_K(S)$ 
is {\em satisfied by} $A$ if for any set map $S\to A$ the induced map
$Mag_K(S)\to A$ sends $\omega$ to $0$. A set of equations
defines a {\em variety} of $K$-algebras, namely the subcategory of $\bMag_K$ 
cut out by (that is, satisfying) these equations. An object of a variety
of $K$-algebras $\bV$ is a ``$\bV$-algebra.'' 

A variety of $K$-algebras is an ``algebraic category'' \cite{frankland}.
It is complete
and cocomplete. Any subalgebra of a $\bV$-algebra is again a $\bV$-algebra.
The forgetful functor $u:\bV\to\bLMod_K$ to the category of left $K$-modules
has a left adjoint
\[
F:\bLMod_K\to\bV\,.
\]

\begin{examples}{\rm Here are four standard examples, beyond $\bMag_K$ itself.
\begin{itemize}
\item $\bAss_K$, the variety of associative algebras, is defined 
by the equation 
\[
(xy)z-x(yz)\in Mag_K\{x,y,z\}\,.
\]
\item Adding the further equation 
\[
xy-yx\in Mag_K\{x,y\}
\]
gives us the variety of commutative $K$-algebras, $\bCom_K$. 
\item A Lie algebra (in $\bLie_K$) 
is a $K$-algebra satisfying the equations
\[
xx\in Mag_K\{x\}\,,\quad (xy)z+(yz)x+(zx)y\in Mag_K\{x,y,z\}\,.
\]
\item An alternative algebra is a magmatic $K$-algebra satisfying the 
equations 
\[
(xx)y-x(xy)\,,\, (xy)y-x(yy)\,\in\,Mag_K\{x,y\}\,.
\]
These are the objects in the variety $\bAlt_K$.

\end{itemize}}
  
\end{examples}

Note that we do not assume a unit element in any of these examples.
The example of alternatiive algebras is less familiar than the others and 
we spend a moment introducing it. 
A great deal of information about alternative algebras can be found in 
the books of Richard Schaefer \cite{schafer} and Kevin McCrimmon 
\cite{mccrimmon}.
Any associative $K$-algebra is alternative. If $A$ is a magmatic $K$-algebra
such that every pair of elements generates an associative subalgebra, then
$A$ is alternative; and Emil Artin proved the converse statement
\cite{bruck-kleinfeld}, \cite[\S3.3]{mccrimmon}.  
The algebra of octonions \cite{baez} is a well-known example of a
nonassociative alternative algebra. 

A map of sets $S\times S\to S$ induces, by linear extension, a magmatic
algebra structure on $KS$. If the product on $S$ is associative, then $KS$
is an associative algebra. But the analogous assertion is false for the
alternative axioms. For example the multiplication on $\{a,b,c\}$ given by
\begin{center}
\begin{tabular}{c|ccc}
& $a$ & $b$ & $c$ \\
\hline
$a$ & $a$ & $a$ & $c$ \\
$b$ & $a$ & $b$ & $b$ \\
$c$ & $c$ & $b$ & $c$ 
\end{tabular}
\end{center}
is commutative and alternative, but the algebra
structure on the $K$-module that it generates is not alternative;
for example $((a+b)(a+b))c\neq(a+b)((a+b)c)$. 
We thank Hadeel AbuTabeekh for this example. 

\begin{remark}{
There is a reversal involution $\overline{(-)}:\bMag\to\bMag$. It comes with
a natural bijection of underlying sets $X\to\overline X$ that we will also
denote with an 
overline, and $\overline x\,\overline y=\overline{yx}$. It extends to an 
involution of $\bMag_K$. To any variety $\bV$ of $K$-algebras we can associate
an ``opposite'' variety $\overline\bV$, with defining equations given by 
reversing the defining equations of $\bV$. By sending a $\bV$-algebra to
the same $K$-module with opposite multiplication, we get a natural 
equivalence of categories
\[
\bV\to\overline\bV\,,\quad A\mapsto \overline A\,.
\]
A variety $\bV$ is {\em symmetric} if
$\overline\bV=\bV$. All the examples above are symmetric,
but, for example, the variety of ``left alternative $K$-algebras,'' satisfying 
$(xx)y-x(xy)$ but perhaps not $(xy)y-x(yy)$, is not symmetric; its
opposite is the variety of right alternative $K$-algebras. 
If $\bV$ is a symmetric variety, the isomorphism $\bV\to\overline\bV$ 
becomes an involution on $\bV$, sending an algebra to the same $K$-module
with the opposite multiplication.
}\end{remark}

\section{Beck modules}
\label{sec-beck}

Let $\bV$ be a variety of $K$-algebras and $A$ a $\bV$-algebra. 
The ``slice category'' $\bV/A$ has as objects morphisms in $\bV$ with
target $A$, and as morphisms maps compatible with the structure maps to $A$. 
This slice category again has good properties; in particular it is complete
and cocomplete. We can thus speak of abelian group objects in $\bV/A$.

An abelian
group structure on a $\bV$-algebra over $A$, $p:B\downarrow A$, begins with
a unit: a map from the terminal object of $\bV/A$, that is, a section 
$\eta:A\uparrow B$ of $p$. This unit defines an ``axis inclusion''
$i:B\coprod B\to B\times_AB$ in $\bV/A$. 
A magma structure on $B$ is an extension
of the ``fold map'' $\nabla:B\coprod B\to B$ over the product. In these
algebraic situations, the map $i$ is an epimorphism, so such an extension
is unique if it exists: Being a unital magma object in $\bV/A$ is a 
{\em property} of a pointed object, not further structure on it.
Furthermore, the unique unital magma structure with given unit, when
it exists, is an abelian group structure. We call an object of this type
an {\em abelian object}.

\begin{definition} \cite{beck}
Let $A$ be a $\bV$-algebra. A {\em Beck $A$-module} 
is an abelian object in the slice category $\bV/A$:
\[
\bMod_A=\Ab(\bV/A)\,.
\] 
\end{definition}

The trivial $K$-module has a unique $\bV$-algebra structure and forms the
terminal object in $\bV$. A Beck $0$-module -- an abelian object in $\bV$ -- 
is a $K$-module equipped with the trivial product, $ab=0$ for all $a,b$. 

\begin{proposition}{\cite[Theorem 3.16]{adamek-rosicky} and
    \cite[Chapter 2, Theorem 2.4]{barr}} $\bMod_A$ is a complete and
cocomplete abelian category.
  \end{proposition}

It's clear that if $f:M\to N$ is a morphism in $\bMod_A$ and $c\in K$ then
the map $cf\in\bMag_K(M,N)$ is again a morphism of $A$-modules: $\bMod_A$ 
is $K$-linear. We will denote the $K$-module of morphisms $M\to N$ in 
$\bMod_A$ by $\Hom_A(M,N)$. 

In our $K$-linear situation, write $M$ for the kernel of $p:B\downarrow A$.
Suppose $B\downarrow A$ has the structure of a unital magma in $\bV/A$.
This consists
of two pieces of structure: the ``unit'' is a map from the terminal object in 
$\bLMod_K/A$, that is, a section of $p:B\downarrow A$, and the ``addition,''
a map $\alpha:B\times_AB\to B$ over $A$. Since 
\[
B\times_AB=(A\oplus M)\times_A(A\oplus M)=A\oplus M\oplus M
\]
the structure map has the form $\alpha:A\oplus M\oplus M\to A\oplus M$.
Using linearity and unitality it's easy to see that the 
``addition'' is actually determined by the addition in $M$:
\[
\alpha(a,x,y)=(a,x+y)\,.
\]
The $K$-algebra structure on $A\oplus M$ is described by left and right 
``actions''
\[
A\tensor M\to M\,,\quad M\tensor A\to M
\]
both of which we denote by juxtaposition. Together they determine the 
multiplication on $A\oplus M$ by
\[
(a,x)(b,y)=(ab,ay+xb)\,.
\]

Absent further axioms, these action maps satisfy no properties. 
This describes the category of magmatic Beck $A$-modules. It is
equivalent to the category of right modules over $\Tens_K(A\oplus A)$.
Let $\lambda:A\to\Tens_K(A\oplus A)$ denote the inclusion of the first
factor, and $\rho$ the inclusion of the right factor. Then the action
of $\Tens_K(A\oplus A)$ on $M$ is given by
\[
x\lambda(a)=ax\,,\quad x\rho(a)=xa\,.
\]

If we are working with a 
general variety of $K$-algebras $\bV$, the axioms of $\bV$ will determine
further properties of these two actions. For example, with $\bV=\bAss$, 
these left and right ``actions'' are required to satisfy 
\[
(xb)c=x(bc)\,,\quad (ay)c=a(yc)\,,\quad (ab)z=a(bz)
\]
for all $a,b,c\in A$ and $x,y,z\in M$. 
In other words, $\bMod_A$ is the usual category of bimodules over $A_+$
(for which $K$ acts the same way on both sides), where 
\[
A_+=K\oplus A
\]
with product given by $(p,a)(q,b)=(pq,pb+qa)$ is the unital $K$-algebra 
associated to $A$.

Forming the underlying $K$-module of a Beck $A$-module gives a functor
\[
u:\bMod_A\to\bLMod_K\,.
\]

\begin{lemma} The functor $u$ has a left adjoint
\[
F_A:\bLMod_K\to\bMod_A
\]
sending a $K$-module $V$ to the ``free $A$-module generated by $V$.''
\end{lemma}
\begin{proof}
We appeal to the Freyd adjoint functor theorem, \cite[p. 117]{maclane}.
The functor $u$ reflects limits: If $M:D\to\bMod_A$ is a diagram of 
$A$-modules,
the limit of $uM:D\to\bLMod_K$ has a unique $A$-module structure that serves
as the limit in $\bMod_A$. The solution set condition is this: For any
$V\in\bLMod_K$, we require a set $\Sigma$ of pairs $(M,f)$, where $M\in\bMod_A$
and $f:V\to uM$, with the property that for any $g:V\to uN$ there is 
$(M,f)\in\Sigma$ and an $A$-module map $t:M\to N$ such that $g=(ut)\circ f$. 
To describe an appropriate set $\Sigma$, we employ the following language.
A $V$-{\em generated} $A$-module is a pair $(M,f)$ where $M$ is an $A$-module,
$f:V\to uM$, and $M$ is the minimal $A$-module containing the image of $f$.
We claim that for a given $K$-module $V$, there is a set of $V$-generated
$A$-modules such that any $V$-generated $A$-module is isomorphic to a member
of this set. This is clear if $\bV=\bMag_K$, since then the category of
$A$-modules is equivalent to the category of right modules over the 
unital associative $K$-algebra $\Tens_K(A\oplus A)$, and a $V$-generated right
module over this $K$-algebra is isomorphic to one of the form 
\[
V\to V\tensor_K\Tens_K(A\oplus A)\to M
\]
where the first map sends $v$ to $v\tensor 1$ and the second is a surjection
of right modules. 
But if $A\in\bV$, then the isomorphism
classes of $V$-generated $A$-modules in $\bV$ are among the isomorphism
classes of $V$-generated $A$-modules in $\bMag_K$, and so also form a set.

We can now take $\Sigma$ to be a set of representatives of isomorphism
classes of $V$-generated Beck $A$-modules.
\end{proof}

Spelling out the adjunction, we have a bijection, natural in the 
pair $V\in\bLMod_K$ and $M\in\bMod_A$:
\[
\Hom_A(F_AV,M)=\Hom_K(V,uM)\,.
\]
In particular, 
\[
\Hom_A(F_AK,M)=uM\,.
\]
Since $u$ is exact, the object $F_AK$ is a projective generator of $\bMod_A$.
This lets us apply another theorem of Freyd's, the embedding theorem
\cite[p. 106]{freyd}, to identify the category $\bMod_A$ with the category
of right modules over a certain unital associative $K$-algebra.

\begin{definition}{\em 
Let $\bV$ be a variety of $K$-algebras and $A\in\bV$. The {\em universal
enveloping algebra} of $A$, $U_\bV(A)$, is the unital associative $K$-algebra
\[
U_\bV(A)=\End_A(F_AK)\,.
\]
}\end{definition}

For any $M\in\bMod_A$, the $K$-module underlying $M$ thus admits a natural
right module structure over $U_\bV(A)$, given by precomposing with the
endomorphism of $F_A(K)$. To summarize:

\begin{proposition} 
This construction provides a natural equivalence of abelian 
categories
\[
\bMod_A\to\bRMod_{U_\bV(A)}\,. \qedhere
\]
\end{proposition}

This construction is of course natural in the $\bV$-algebra $A$. Any variety of
$K$-algebras has a terminal object, the trivial $K$-module $0$, and the
defining property of the universal enveloping algebra implies that
$U_\bV(0)=K$. So a universal enveloping algebra always has a canonical 
augmentation 
\[
\epsilon:U_\bV(A)\to K\,.
\] 

Let $\Omega'$ and $\Omega$ be two sets of equations, cutting out varieties
$\bV'$ and $\bV$ of $K$-algebras. If $\Omega'\subseteq\Omega$, then any
$\bV$-algebra is a $\bV'$-algebra; write $i:\bV\to\bV'$ for
the inclusion functor. Fix $A\in\bV$. There is then a functor
\[
i_*:\bMod_{A}\to\bMod_{iA}
\]
that sends $M$ to itself as a $K$-module, with the same left and right actions
by $A$ but now regarded as giving a Beck $iA$-module structure. This 
functor is induced by a map of unital associative algebras
\[
i^*:U_{\bV'}(iA)\to U_\bV(A)\,.
\]

In particular we might take $\Omega'$ to be empty, so that $\bV'=\bMag_K$
and $U_{\bV'}(A)=\Tens_K(A\oplus A)$. For any $\bV$ we thus receive 
a canonical map of unital associative $K$-algebras
\[
\pi:\Tens_K(A\oplus A)\to U_\bV(A)\,.
\]
Denote the composite $\pi\circ\lambda$ by $l:A\to U_\bV(A)$ and
$\pi\circ\rho$ by $r:A\to U_\bV(A)$. These are $K$-linear maps, and the
sum of their images generates $U_\bV(A)$ as an associative unital $K$-algebra;
the map $\pi$ is surjective.

\section{Noncommutative partial differentiation}
\label{sec-partial-diff}

Let $\bV$ be a variety of $K$-algebras.
The universal enveloping algebra of a $\bV$-algebra $A$ is the quotient of
$U_{\bMag_K}(A)=\Tens_K(A\oplus A)$ by an ideal generated by elements 
determined
by the equations defining $\bV$. These elements are derived by a process of 
``noncommutative partial differentiation,'' which we now describe.

Let $Mag(S)$ be the free magma on a set $S$, and $Ass(X)$ the free
associative algebra on the set $X$. Write
\[
\lambda\,,\,\rho:Mag(S)\to Ass(Mag(S)\sqcup Mag(S))
\]
for the inclusions of the left and right summands. 

\begin{lemma}
For each $x\in S$ there is a unique map
\[
\frac{\partial}{\partial x}:Mag(S)\to Ass(Mag(S)\sqcup Mag(S))
\]
such that
\[
\frac{\partial x}{\partial x}=1\,;
\]
for any $y\in S$ with $y\neq x$
\[
\frac{\partial y}{\partial x}=0\,;
\]
and for any $\alpha,\beta\in Mag(S)$,
\[
\frac{\partial\alpha\beta}{\partial x}=
  \frac{\partial\alpha}{\partial x}\rho_\beta+
  \frac{\partial\beta}{\partial x}\lambda_\alpha\,.
\]
\end{lemma}

\begin{proof} This is immediate, since any element of $Mag(S)$ is 
built up from elements of $S$ by a unique sequence of multiplications. 
\end{proof}

For any commutative ring $K$, this map extends by linearity to
\[
\frac{\partial}{\partial x}:KMag(S)\to\Tens_K(KMag(S)\oplus KMag(S))\,.
\]

\begin{example}{\em 
Define the left- and right-bracketed powers of $x$ by
\[
x^{(1}=x\,,\, x^{(n}=x^{(n-1}x\,,\quad x^{1)}=x\,,\, x^{n)}=xx^{n-1)}
\]
Then for $n>1$ we have the following noncommutative analogues of the familiar
formula for the derivative of a power:
\begin{gather*}
\frac{\partial x^{(n}}{\partial x}=\lambda_{x^{(n-1}}+\lambda_{x^{(n-2}}\rho_x
    + \cdots + \lambda_x\rho_x^{n-2}+\rho_x^{n-1} \\
    \frac{\partial x^{n)}}{\partial x}=\rho_{x^{n-1)}}+\rho_{x^{n-2)}}\lambda_x
    + \cdots + \rho_x\lambda_x^{n-2}+\lambda_x^{n-1} \,.
\end{gather*}
}\end{example}

The application of this operation to the determination of the 
structure of universal enveloping algebras is this: 
Recall that a variety of $K$-algebras is cut out by a set $\Omega$ of 
equations. An element of $\Omega$ is a pair $(S,\omega)$ where $S$ is 
a finite set and $\omega\in KMag(S)$. If $A$ is a $K$-algebra, a set map 
$a:S\to A$ determines a map of $K$-algebras $KMag(S)\to A$, which we denote
by $\omega\mapsto\omega(a)$. Putting this map on both factors gives us
$KMag(S)\oplus KMag(S)\to A\oplus A$, and hence to a map of associative
unital $K$-algebras
\[
\Tens_K(KMag(S)\oplus KMag(S))\to\Tens_K(A\oplus A)
\]
which we denote again by $\mu\mapsto\mu(a)$. Here is the description 
of the ideal defining $U_\bV(A)$.

\begin{proposition} 
Suppose $\bV$ is a variety of $K$-algebras cut out by a set of equations 
$\Omega$. Let $A$ be a $\bV$-algebra. Then the universal enveloping 
algebra $U_\bV(A)$ is the quotient of $\Tens_K(A\oplus A)$ by the ideal $I$
generated by the set 
  \[
\left\{\frac{\partial\omega}{\partial x}(a)\,:\,
  (S,\omega)\in\Omega\,,\, x\in S\,,\, a:S\to A\right\}\,.
  \]
\end{proposition}

\begin{proof} Let $A$ be a magmatic $K$-algebra and $M$ a module for it:
so $A\oplus M$ has the structure of a magmatic $K$-algebra with product
given by $(a,m)(b,n)=(ab,an+mb)$. $M$ is then a right module for the 
magmatic universal enveloping algebra
\[
U_\bMag(A)=\Tens_K(A\oplus A)\,.
\]
Given $\omega\in KMag(S)$, we claim that for any $(a,m):S\to A\oplus M$: 
\[
\omega(a,m)=
\left(\omega(a),\sum_{x\in S}m_x\frac{\partial\omega}{\partial x}(a)\right)
\in A\oplus M\,.
\]

Since the definition of the partial derivative is inductive, we proceed by
induction. For the base case, we note that it is true if $\omega=x$ for 
some $x\in S$:
then $x(a,m)=(a_x,m_x)$, which agrees with the right hand side since 
all but one term vanishes in the sum. Then, given $\alpha,\beta\in K\bMag(S)$, 
we compute
\begin{align*}
\alpha\beta(a,m)= & \,\, \alpha(a,m)\beta(a,m)\\
= & \left(\alpha(a),\sum m_x\frac{\partial\alpha}{\partial x}(a)\right)
\left(\beta(a),\sum m_x\frac{\partial\beta}{\partial x}(a)\right) \\
= & \left(\alpha(a)\beta(a),
\sum m_x\frac{\partial\beta}{\partial x}\lambda_{\alpha(a)}
+\sum m_x\frac{\partial\alpha}{\partial x}\rho_{\beta(a)}\right) \\
= & \left(\alpha\beta(a),
\sum m_x\left(\frac{\partial\beta}{\partial x}(a)\lambda_{\alpha(a)}
+\frac{\partial\alpha}{\partial x}(a)\rho_{\beta(a)}\right)\right)
\end{align*}
and the factor in the sum is indeed 
$\displaystyle{\frac{\partial\alpha\beta}{\partial x}(a)}$.

Now suppose that $\bV$ is cut out by $\Omega$, that $A$ is a $\bV$-algebra,
and that $M$ is a Beck $A$-module, and let $\omega\in\Omega$ and 
$(a,m):S(\omega)\to A\oplus M$. Since $A\oplus M$ is a $\bV$-algebra,
$\omega(a,m)=0$, so the right hand side of the above equation vanishes.
Let $x\in S$, and take $m:S\to M$ be a function that vanishes except at $x$.
We discover that $\partial\omega/\partial x(a)$ vanishes on $M$. Since
$M$ was an arbitrary Beck $A$-module, this element lies in the ideal $I$. 

Since the noncommutative partial derivatives of the defining equations for
$\bV$, evaluated on maps to $A$, simply record the validity of those 
equations on a Beck $A$-module, there are no further relations in the 
ideal $I$.
\end{proof}

\begin{example}{\em The equation $xy$ cuts out the variety of $K$-algebras with
trivial multiplication; this is just the category $\bLMod_K$ of $K$-modules.
We find 
\[
\frac{\partial xy}{\partial x}=\rho_y\quad,\quad\quad
\frac{\partial xy}{\partial y}=\lambda_x
\]
so for any algebra $A$ in this variety, $U(A)=K$; the category of Beck
$A$-modules is again just the category of $K$-modules.
}\end{example}

\begin{example}{\em  $\bAss_K$ is cut out by $(xy)z-x(yz)$. Compute:
  \begin{gather*}
  \frac{\partial(xy)z}{\partial x}=\frac{\partial xy}{\partial x}\rho_z
  =\rho_y\rho_z \quad,\quad\quad
  \frac{\partial x(yz)}{\partial x}=\rho_{yz}\\
  \frac{\partial(xy)z}{\partial y}=\frac{\partial xy}{\partial y}\rho_z
  =\lambda_x\rho_z\quad,\quad\quad
  \frac{\partial x(yz)}{\partial y}=\frac{\partial yz}{\partial y}\lambda_x
  =\rho_z\lambda_x\\
  \frac{\partial(xy)z}{\partial z}=\lambda_{xy}\quad,\quad\quad
  \frac{\partial x(yz)}{\partial z}=\frac{\partial yz}{\partial z}\lambda_x
  =\lambda_y\lambda_x
  \end{gather*}
  so the defining relations for $U_\bAss(A)$ are
  \[
  r_ar_b=r_{ab}\,,\quad l_ar_b=r_bl_a\,,\quad l_{ab}=l_bl_a
  \]
  for $a,b\in A$. This shows that Beck $A$-modules in $\bAss_K$ are precisely
$A_+$-bimodules.
The quotient of $\Tens_K(A\oplus A)$ by these relations is the ``extended 
$K$-algebra''
 \[
U_\bAss(A)=A^{op}_+\tensor_KA_+\,.
\]
}\end{example}

\begin{example}{\em Commutativity is specified by $xy-yx$. Compute
    \begin{gather*}
      \frac{\partial xy}{\partial x}=\rho_y\quad,\quad\quad
      \frac{\partial yx}{\partial x}=\lambda_y\\
      \frac{\partial xy}{\partial y}=\lambda_x\quad,\quad\quad
      \frac{\partial yx}{\partial y}=\rho_x
    \end{gather*}
    and both equations give us $l_a=r_a$: The corresponding universal
    enveloping algebra is just $\Tens_K(A)$, independent of the algebra
    structure on $A$. A Beck module for a commutative
    magmatic $K$-algebra is simply a $K$-module $V$ together with a $K$-linear
    map $V\tensor_KA\to V$.

If we combine these two, we get the category $\bCom_K$ of 
commutative, associative, nonunital $K$-algebras. Combining $l_a=r_a$
with the the relations in $U_\bAss(A)$, we find that a Beck $A$-module is
simply an $A_+$-module in the usual sense, and
\[
U_\bCom(A)=A_+\,.
\]
}\end{example}

\begin{example}{\em
    The variety $\bLeib_K$ of Leibniz algebras \cite{loday-pirashvili}
over $K$ is cut out by
    \[
   (xy)z-x(yz)-(xz)y\,.
    \]
Compute:
    \begin{gather*}
    \frac{\partial((xy)z-x(yz)-(xz)y)}{\partial x}=
    \rho_y\rho_z-\rho_{yz}-\rho_z\rho_y \\
    \frac{\partial((xy)z-x(yz)-(xz)y)}{\partial y}=
    \lambda_x\rho_z-\rho_z\lambda_x-\lambda_{xz} \\
    \frac{\partial((xy)z-x(yz)-(xz)y)}{\partial z}=
    \lambda_{xy}-\lambda_y\lambda_x-\lambda_x\rho_y\,.
    \end{gather*}
This leads to the relations in 
$U_{\bLeib}(A)$ \cite[\S2]{loday-pirashvili}:
    \[
    r_{ab}=[r_a,r_b]\,,\quad l_{ab}=[l_a,r_b]\,,\quad (l_a+r_a)l_b=0
    \]
(where the bracket denotes the commutator in this associative algebra). 
    If we adjoin the relation $xx$, to get the variety $\bLie_K$,
    we find $\lambda_x=-\rho_x$
    and so $l_a=-r_a$ in $U_{\bLie}(A)$, and $U_\bLie(A)$ is the quotient of
    $\Tens_K(A)$ by the relations
    \[
    r_{ab}=[r_a,r_b]\,,
    \]
giving us the usual Lie universal enveloping algebra.
}\end{example}

\begin{example}{\em
Chataur and Livernet \cite{chataur-livernet} 
consider ``level algebras,'' defined by the equations
\[
xy-yx\,,\, (wx)(yz)-(wy)(xz)\,.
\]
The second equation leads to
\[
\rho_x\rho_{yz}=\rho_y\rho_{xz}\,,\,
\lambda_w\rho_{yz}=\rho_z\lambda_{wy}\,,\,
\rho_z\lambda_{wx}=\lambda_w\rho_{xz}\,,\,
\lambda_y\lambda_{wx}=\lambda_x\lambda_{wy}\,.
\]
Commutativity leads to $\lambda_x=\rho_x$, and hence to the attractive 
relations in $U_\bLev(A)$:
\[
U_\bLev(A)=\Tens_K(A)/(r_ar_{bc}=r_br_{ca}=r_cr_{ab}\,,\, a,b,c \in A)\,.
\]
Since they are multilinear, the relations can be restricted to hold for
$a,b,c$ belonging to a basis for $A$. The universal enveloping algebra
of an abelian level algebra is the tensor algebra on the underlying module.
In general, the relations are homogeneous; the universal algebra is always a 
quadratic graded algebra. We will return to this interesting example in
later work.
}\end{example}

\begin{example} {\em 
Even the variety of $K$-algebras cut out by $(xx)x$ has 
interesting universal enveloping algebras: The relations are
\[
r_a^2+l_ar_a+l_{aa}=0\,.
\]
Replacing $a$ by $a+b$ gives
\[
r_ar_b+r_br_a+l_ar_b+l_br_a+l_{ab}+l_{ba}=0\,,
\]
and since these relations are multilinear it suffices to require them
for $a,b$ in a set of $K$-module generators for $A$.
In the abelian case of $K^n$, with basis $e_1,\ldots,e_n$, let 
$l_i=l_{e_i}$ and $r_i=r_{e_i}$. The relations are neater if we let 
$q_i=l_i+r_i$, for then they are
\[
q_ir_i\,\,\text{for}\,\,1\leq i\leq n\,\,,\quad 
q_ir_j+q_jr_i\,\,\text{for}\,\, 1\leq i<j\leq n\,.
\]
}
\end{example}

\begin{example} \label{alternative-relations}
{\em The alternative equations differentiate to
\begin{gather*}
\rho_{xy}-\rho_x\rho_y=\lambda_x\rho_y-\rho_y\lambda_x\,,\quad 
\lambda_{xx}=\lambda_x\lambda_x\,,\\
\rho_y\rho_y=\rho_{yy}\,,\quad
\lambda_{xy}-\lambda_y\lambda_x=\rho_y\lambda_x-\lambda_x\rho_y\,,
\end{gather*}
so $U_{\bAlt}(A)$ is $\Tens_K(A\oplus A)$ modulo the ideal generated by
\begin{gather*}
r_{bb}=r_br_b\,,\quad l_{aa}=l_al_a\,,\\
l_{ab}-l_bl_a=r_bl_a-l_ar_b=r_ar_b-r_{ab}\,.
\end{gather*}
}\end{example}

\begin{remark}{\em 
This use of the term ``universal enveloping algebra'' differs from the 
classical Lie perspective. In that case, one has a ``forgetful'' functor
$\bAss_K\to\bLie_K$ that sends an associative $K$-algebra to the Lie algebra
structure on the $K$-module $A$ given by $[a,b]=ab-ba$; and $U$ is the left
adjoint of this functor. In our generality, there is no ``underlying''
associative algebra; the universal enveloping algebra has a different 
defining property. But it turns out to produce the same result in the case
of $\bLie_K$. Our meaning for the term was explored in the operadic case
by Ginzburg and Kapranov in \cite{ginzburg-kapranov}, and in the example
of Leibniz algebras by Loday and Pirashvili in \cite{loday-pirashvili}.
For alternative algebras, the universal enveloping algebra was implicit 
in the work of Jacobson, Stanley, 
and others, and explicit in Kevin McCrimmon's unpublished book 
\cite[1.7, p.~7-9]{mccrimmon}, where it is called the ``universal 
multiplication algebra'' of $A$. His comment about this construction:
``In practice we will try to avoid the universal multiplication algebra
and deal as much as possible with ``real'' multiplications, because it is
more comfortable dealing with honest-to-goodness multiplication operators
than than with the highly formal ``universal'' or ``ideal'' multiplications.''
This use of the term has become standard in homotopy theory; see
\cite{knudsen} for example.
}\end{remark}
  
\begin{remark} \label{rem-compare}
{\em
We have three closely related varieties of $K$-algebras, related by
forgetful right adjoints
\[
\bCom_K\ra{u}\bAss_K\ra{u}\bAlt_K\,.
\]
The left adjoints are given by forming the maximal associative quotient of an
alternative algebra, and the maximal commutative quotient of an associative
algebra. These functors induce right adjoints
\[
\bCom_K/A\to\bAss_K/uA
\]
for $A\in\bCom_K$ and
\[
\bAss_K/A\to\bAlt_K/uA
\]
for $A\in\bAss_K$. As right adjoints, they preserve products, and hence induce
functors
\[
\Ab(\bCom_K/A)\to\Ab(\bAss_K/uA)\,,\quad\Ab(\bAss_K/A)\to\Ab(\bAlt_K/uA)\,.
\]

These functors can be described by means of ring homomorphisms between the 
corresponding universal enveloping algebras:
There is a unital $K$-algebra surjection natural in $A\in\bAss_K$ 
\[
U_{\bAlt}(uA)\to U_{\bAss}(A)\,,\quad 
l_a\mapsto a\tensor 1\,,\quad r_b\mapsto1\tensor b\,,
\]
and a unital $K$-algebra surjection natural in $A\in\bCom_K$ 
\[
U_{\bAss}(uA)\to U_{\bCom}(A)\,,\quad a\tensor b\mapsto ab\,.
\]
}\end{remark}

\begin{remark}{\em 
The reversal endomorphism of $\bMag$ induces natural isomorphisms
\[
U_{\overline\bV}(\overline A)=U_{\bV}(A)^{op}
\]
that swaps the $K$-module maps $r,l$ from the $K$-module $A$. 
}\end{remark}

Formation of the universal enveloping algebra enjoys a strong base-change
property. 
A homomorphism of commutative rings $f:K\to L$ induces a multiplicative map
$f:KMag(S)\to LMag(S)$. By applying this map to the equations defining
a variety $\bV_K$ of $K$-algebras, we obtain a base-changed variety
$\bV_L$ of $L$-algebras, and a functor $L\tensor_K-:\bV_K\to\bV_L$.

\begin{proposition}\label{prop-basechange}
Fix this notation, and let $A\in\bV_K$. There is a natural map
  $U_{\bV_K}(A)\to U_{\bV_L}(L\tensor_KA)$ of unital $K$-algebras
  that extends to an isomorphism
  \[
L\tensor_KU_{\bV_K}(A)\to U_{\bV_L}(L\tensor_KA)
\]
of unital $L$-algebras.
\end{proposition}

\begin{proof}
To begin with, $f:K\to L$ induces a multiplicative map 
\[
f_*:\Tens_K(KMag(S)\oplus KMag(S))\to\Tens_L(LMag(S)\oplus LMag(S))\,.
\]
Base-change is compatible with partial differentiation: 
Given $\omega\in K\bMag(S)$ and $x\in S$,
\[
f_*\frac{\partial\omega}{\partial x}=\frac{\partial f\omega}{\partial x}\,.
\]
So for any $\bV_K$-algebra $A$, $f_*$ carries the equations defining 
$U_{\bV_K}(A)$ to the equations defining $U_{\bV_L}(L\tensor_KA)$.
We obtain the desired homomorphism on universal enveloping algebras.
Moreover, the relations in $\Tens_K(A\oplus A)$ defining $U_{\bV_K}(A)$
are carried precisely to the relations in 
$\Tens_L((L\tensor_KA)\oplus(L\tensor_KA))$ defining $U_{\bV_L}(L\tensor A)$,
so the result follows.
\end{proof}

\section{Universal enveloping algebras for abelian alternative algebras}
\label{sec-aa-univ}

Let $K^n$ denote the free $K$-module of rank $n$ regarded as an alternative
$K$-algebra with trivial product. 
Write $l_i$ and $r_i$, $1\leq i\leq n$, for the images in $U(K^n)$ of the
standard basis elements under $l,r:K^n\to U(K^n)$. 
This associative $K$-algebra is graded by weight. Clearly $U(K^0)=K$ and
$U(K^1)$ has basis $\{1,l_1,r_1,l_1r_1\}$; in fact
$l_1^2=r_1^2=l_1r_1-r_1l_1=0$.

\begin{theorem} Let $K$ be any commutative ring. For each $n$ and $k$,
the $K$-module $U(K^n)_k$ is free, and
\[
\rank_K U(K^n)_k=
\begin{cases} 1 & \text{if}\,\,k=0\\
2n & \text{if}\,\, k=1\\
\frac{3n^2-n}{2} & \text{if}\,\, k=2\\
2\textstyle{\binom{n}{k}} & \text{if}\,\,k\geq3\,.
\end{cases}
\]
\end{theorem}

In particular, $U(K^n)_k=0$ for $k>n$ as long as
$n>1$. The growth of $\dim_KU(K^n)$ with $n$ is exponential;
\[
\dim_KU(K^n)=2\cdot2^n+\frac{n^2+n-2}{2}\,.
\]

\begin{proof}
In fact we can make a more precise statement, specifying for each $k$ a set 
of monomials in  $S=\{l_1,r_1,\ldots,l_n,r_n\}$ forming a basis for $U(K^n)_k$.
The first step is to determine a Gr\"obner basis for the ideal
\[
I=\ker(\Tens_K(K^n\oplus K^n)\to U(K^n))\,.
\] 
Order $S$ as shown, so that $l_1>r_1>l_2>\cdots$; 
order monomials first by weight and within a given
weight left-lexicographically. The ``leading monomial'' in a polynomial
will be the greatest term. The expressions for the generators of $I$, 
when applied to the elements of $S$, are
\begin{gather*}
l_il_i \,\text{and}\,\, r_ir_i \,\text{for all}\,\, i \\
l_ir_i-r_il_i \,\,\text{for all} \,\, i \\
l_jl_i+r_ir_j \,\,\text{and}\,\, r_jr_i+l_il_j\,\,\text{for}\,\, j<i \\ 
l_jr_i+r_jr_i-r_il_j\,\,\text{and}\,\,r_jl_i-l_ir_j-r_ir_j\,\,\text{for}\,\,j<i
\end{gather*}
The first identities do not imply that $l_a^2=0=r_a^2$ for all $a$, however, 
and to get a complete set of relations involving just the basis elements 
we have to adjoin
$l_jl_i+l_il_j$ and $r_jr_i+r_ir_j$ for $j<i$. Comparing leading terms with 
the existing relations, we can drop $l_jl_i+l_il_j$ and $r_jr_i+l_il_j$ for
$j<i$ at the expense of adjoining $l_il_j-r_ir_j$ for $j<i$. 
We can also use $r_jr_i+r_ir_j$ to replace $l_jr_i+r_jr_i-r_il_j$ 
with $l_jr_i-r_ir_j-r_il_j$, so that all the relations except one have
the effect of replacing an increasing sequence of subscripts by a decreasing
one. We then have a ``reduced'' set of relations,
in the sense that there are no repeated leading terms:
\begin{itemize}
\item weight 2, length 1: $l_il_i$ and $r_ir_i$ for all $i$
\item weight 2, length 2: $l_ir_i-r_il_i$ for all $i$
\item weight 2, length 2: $l_jl_i+r_ir_j$, $r_jr_i+r_ir_j$, 
and $l_il_j-r_ir_j$ for $j<i$ 
\item weight 2, length 3: $l_jr_i-r_il_j-r_ir_j$ and $r_jl_i-l_ir_j-r_ir_j$ 
  for $j<i$
\end{itemize}

Write $R_0$ for this basis for $I$. Notice that there are monomials in the
ideal that are not divisible by any leading entry in this list: for example
\[
r_ir_jl_j=l_i(l_jl_j)-(l_il_j-r_ir_j)l_j\,.
\]
So this is not a Gr\"obner basis. To get one, we need to
inspect overlaps and adjoin any overlap differences that can't be reduced
to zero using the leading monomials of elements of $R_0$, and then repeat that
process if necessary. This is not hard to do by hand, and it leads to the
following additional basis elements. 
\begin{itemize}
\item weight 3, length 1: $r_ir_jl_j,\, r_il_ir_j,\, l_ir_jl_j$ \,for $j<i$ 
\item weight 3, length 2: 
$l_ir_jl_k-r_il_jr_k,\, l_ir_jr_k-r_ir_jl_k$ \,for $k<j<i$
\item weight 3, length 3: $r_il_jr_k+r_ir_jl_k+r_ir_jr_k$ \, for $k<j<i$
\end{itemize}

Let $R_1$ be the set $R_0$ with these new weight 3 relations adjoined. 
One must now inspect overlaps of the new relations with the old ones and with
each other. The result of a hand calculation is that all the overlap 
differences reduce to 0 using $R_1$; we have obtained a Gr\"obner basis 
for $I$. Both these calculations were checked using 
a Python script for $n\leq 5$. This suffices to verify the general case, since
at most 5 indices are involved in any overlap of elements of $R_1$, 
and the $n=5$ case is general enough to cover all possibilities. 

Suppose now that $K$ is a field.
As described in \S\ref{groebner}, the set of ``normal monomials,'' that is, 
those not divisible
by the leading entry of any element of the Gr\"obner basis, projects to 
a vector space basis for the quotient $K$-algebra 
$U(K^n)=\Tens_K(K^n\oplus K^n)/I$. 

Since all the relations are of weight at least 2, we find that $\{1\}$ 
is a basis for $U(K^n)_0$ and $\{l_1,r_1,\ldots,l_n,r_n\}$ is a basis
for $U(K^n)_1$.

The first relation forbids repeated letters. The second implies that if a 
subscript is repeated it must be in the order $rl$. The third and fourth 
relations force the indices in a normal monomial to be decreasing, and
the $l_il_j$ relations imply that there can be no repeated $l$'s.
The other relations are of weight 3, so we find that $U(K^n)_2$ has basis
\[
\{r_il_i\}\sqcup\{l_ir_j,\,r_il_j,\,r_ir_j\,:\,i<j\}
\]
Thus 
\[
\dim_KU(K^n)_2=n+3\binom{n}{2}=\frac{3n^2-n}{2}\,.
\]

The first and third weight 3 length 1 relations imply that  $r_jl_j$ can only 
occur at the beginning of a normal monomial, but since $ll$ can never
occur the second relation in that list rules that out as well. So all
normal monomials have strictly increasing subscripts.

Suppose an $l$ occurs in a monomial of weight at least 3 but not at the end.
If it is the next-to-last entry, it must be in the pattern $rlr$. 
If it is earlier, it must be in $rlr$ or $lrx$ for $x$ either $r$ or $l$. 
All these are excluded by the final three relations. 

So we discover that the normal monomials of weight $k\geq 3$ are exactly 
those with strictly increasing subscripts and which consist entirely of 
$r$'s except possibly for the last entry. There are $2\binom{n}{k}$ of these.

This concludes the proof in case $K$ is a field. For the general case, 
note first that since the tensor algebra is a finitely generated $K$-module
in each weight, we know that $U(K^n)$ is too. Suppose now that $K=\ZZ$. 
By base change \ref{prop-basechange}, $U(\ZZ^n)\tensor\QQ=U(\QQ^n)$, and 
for any prime $p$, $U(\ZZ^n)\tensor\FF_p=U(\FF_p^n)$.
Consequently, if $U(\ZZ^n)_k$ had a nonzero element of order $p$, the 
ranks of $U(\FF_p^n)_k$ and $U(\QQ^n)_k$ would differ. But our calculation
showed that these dimensions were independent of the field. 

Finally, we can appeal to base-change again to pass to an arbitrary
commutative base ring $K$. 
\end{proof}

We end by pointing out that this theorem provides, for any natural number $n$,
an example of an alternative $K$-algebra $E_n$, free of finite rank as a 
$K$-module, with the property that $(uv)(xy)=0$ for all $u,v,x,y\in E_n$ but
containing elements $x_1,x_2,\ldots,x_n$ such that 
\[
(\cdots((x_1x_2)x_3)\cdots x_{n-1})x_n\neq0\,.
\]
In fact we can take for $E_n$ the abelian object in $\bAlt_K/K^n$ 
corresponding to $U(K^n)$ (regarded as a free right module of rank 1 
over itself). So $E_n=K^n\oplus U(K^n)$, with product given by
\[
(a,\alpha)(b,\beta)=(0,\alpha r_b+\beta l_a)\,.
\]
Denote the standard basis for $K^n$ by $\{e_1,\ldots,e_n\}$. Then
$(0,r_1)(e_2,0)=(0,r_1r_2)$, and inductively we find that the indicated
product is nonzero, with 
$x_1=(0,r_1)$ and $x_i=(e_i,0)$ for $i>1$. This is in fact precisely 
the example discovered by G. V. Dorofeev \cite{dorofeev},
\cite[\S6.2]{zhevlakov-et-al}.
That work does not illuminate the role of the universal enveloping algebra,
however.

\section{Universal enveloping algebras for some other alternative algebras}

For any $K$-module $V$, filter the tensor algebra $T=\Tens_K(V)$ by weight. 
Write $d:T-\{0\}\to\NN$ for the function assigning to a nonzero element 
$\alpha$ in the tensor algebra the smallest integer $i$ such that 
$\alpha\in F_iT$. A subset $S\subseteq T-\{0\}$  
can be said to ``satisfy the Poincar\'e-Birkhoff-Witt property'' 
\cite{shepler-witherspoon} if the set
\[
\{s+F_{d(s)-1}T\in\gr_{d(s)}T:s\in S\}
\]
generates the kernel of the induced map 
\[
\gr_*T\to\gr_*(T/I)
\]
where $I$ is the ideal generated by $S$. 

Generally these elements in the associated graded algebra do not generate 
the kernel; but (cf. \cite[Thm. 4.2.1, p. 129]{li}):

\begin{lemma} Let $K$ be a field. 
Any Gr\"obner basis \ref{def-grobner} satisfies PBW.
\end{lemma}

The classical PBW theorem asserts that the image of a basis for a Lie algebra
$A$ over a field $K$
in its universal enveloping algebra satisfies PBW. The consequence is that 
\[
\gr_*U_\bLie(A)=U_\bLie(A^\natural)
\]
where $A^\natural$ is the $K$-module $A$ with trivial Lie bracket. 

For any variety $\bV$ of $K$-algebras and any $\bV$-algebra $A$, the universal
enveloping algebra $U_\bV(A)$ has a canonical filtration by weight: the image
of the weight filtration in $\Tens_K(A\oplus A)$.

In case $\bV=\bAlt_K$, the image of the relations 
\ref{alternative-relations} in the associated graded algebra
$\gr_*U_\bAlt(A)$ are precisely the relations defining $U_\bAlt(A^\natural)$, 
so there is a surjection
\[
U_\bAlt(A^\natural)\to\gr_*U_\bAlt(A)\,.
\]

\begin{corollary} If $A$ is an alternative $K$-algebra that is 
free with $n$ generators as a $K$-module, then
$U_\bAlt(A)$ admits a set of $K$-module generators with at most
\[
2\cdot2^n+\frac{n^2+n-2}{2}
\]
elements.
\end{corollary}

But this upper bound is often not met; this version of the ``PBW theorem'' 
fails in the variety of alternative algebras.

\begin{example} {\em 
Let $A=K[i]/(i^2+1)$, so if $K=\RR$, $A=\CC$. Then the natural
map \ref{rem-compare}
\[
U_\bAlt(A)\to U_\bAss(A)=A_+^{op}\tensor_K A_+
\]
turns out to be an isomorphism. This is easily verified: Inspection of the 
relations defining $U_\bAlt(A)$ shows that it is  generated by $r_i$ and 
$l_i$, that they commute, that each one satisfies $s^3+s=0$, 
and there are no other relations. On the other hand, $A_+=K[s]/(s^3+s)$,
so $U_\bAlt(A)$ is the tensor product of two copies of $A_+$. 
Thus $U_\bAlt(A)$ is free of rank 9, while 
$U_\bAlt(A^\natural)$ is free of rank 10. 
}\end{example}

\begin{example} {\em Computer calculations \cite{yin} with the Buchberger
algorithm lead to the following results. 

Let $A$ be the associative $K$-algebra with generators $i,j$ satisfying
$i^2=j^2=-1,ij=-ji$: the standard quaternion algebra over $K$. Then
$U_\bAlt(A)$ is free of rank 29 as a $K$-module.

Perhaps the most well-known examples of alternative algebras that are not
associative are the octonion algebras \cite{baez}.
Let $A$ be the standard octonion algebra \cite{wiki-octonion}
over $K$. If 2 is invertible in $K$,
$U_\bAlt(A)$ is free of rank 65 as a $K$-module, while if $K$ is a field
of characteristic 2, $\dim_KU_\bAlt(A)=113$.
}\end{example}

\section{Appendix: Gr\"obner bases}
\label{groebner}
We elaborate briefly on the ``Buchberger algorithm'' for expanding a basis
of an ideal in a tensor algebra to a Gr\"obner basis.

Let $S$ be a set equipped with a well-founded partial order: a partial order
in which every strictly decreasing sequence is finite.
The free monoid $B$ generated by $S$ inherits a partial order --
first by weight, and left-lexicographically within a given
weight -- that is again well-founded. 

Let $K$ be a field. The free $K$-module generated by $B$ is the tensor
algebra on $S$, $T=KB$.

Any nonzero element in $T$ has a ``leading monomial,'' the greatest monomial 
occurring with nonzero coefficient in 
its expression as a linear combination of elements of $B$. Write 
\[
LM:T^*\to B
\]
for this function, where we write $I^*=I-\{0\}$ for any ideal $I$ in $T$.
The partial order on $B$ pulls back to a well-founded pre-ordering on $T^*$.

For $u,v\in B$, say that $u$ divides $v$, $u|v$, if there are monomials 
$s,t$ such that $v=sut$. Divisiblity is transitive. 

\begin{definition} \label{def-grobner}
Let $I$ be an ideal in $T$. A subset $G$ of $I^*$ is a 
{\em Gr\"obner basis} for $I$ if for all $r\in I^*$ 
there is $g\in G$ such that $LM(g)|LM(r)$. 
\end{definition}

A Gr\"obner basis $G$ for $I$ yields an efficient algorithm for deciding
whether $z\in T^*$ lies in $I^*$. We may assume $z$ is monic;
that is, the coefficient of $LM(z)$ is $1$.
If $LM(z)$ is not divisible by $LM(g)$ for
any $g\in G$, then $z\not\in I^*$. If instead
\[
LM(z)=sLM(g)t
\]
for some $g\in G$ and $s,t\in B$, then form
\[
z'=z-sgt\,.
\]
This is in $I$ if and only if $z$ is.
If $z'=0$, we have established that $z\in I$.
If not, at least we can say that the leading monomials cancel, so $z'$ is
strictly less than $z$ in the well order on $T^*$.
Divide $z'$ by the coefficient of its leading monomial,
and repeat this process, which terminates because the
ordering is well-founded. The original element $z$ is in $I$ if and only
if the element you wind up with is $0$, in which case you have written $z$
as an explicit linear combination of terms divisible by elements of $G$.
This shows that $G$ generates $I$ as an ideal.

Given a subset $R$ of $I^*$ that generates $I$ as an ideal,
we attempt to enlarge it to a
Gr\"obner basis for $I$. We may assume that each element of $R$ is
monic. Say that  $r,s\in R$ {\em overlap} if there are
$a,b,c\in B$ such that $LM(r)=ab$ and $LM(s)=bc$. For each overlapping
pair $r,s$, form the ``overlap difference''
\[
rc-as\,.
\]
This difference again lies in $I$, but exhibits a new leading monomial,
one that was hidden in the original generating set $R$. 

Now use the reduction process with respect to $R$ to simplify each of the
overlap differences.

Then adjoin all the nonzero reduced overlap differences, made monic,  
to the set $R$ to get a new generating set $R'$. One may want to precede this
step by doing some linear algebra to find a simpler basis for the space
spanned by these polynomials.

\begin{proposition}[\cite{li}, Prop. 5.2, p.~95] 
If $R$ is finite, this process terminates, and the result is a Gr\"obner
basis for $I$.
\end{proposition}

A Gr\"obner basis is minimal (no subset is a Gr\"obner basis for the same
ideal) if and only if it is reduced (no divisibility relations among its
leading monomials) \cite[p.~67]{li}. One may always refine a Gr\"obner
basis to a reduced one. 

A monomial $u$ is {\em normal} mod $I$ if it is not divisible by any element
of $LM(I^*)$. If $G$ is a Gr\"obner basis for $I$, it suffices to check
non-divisibility by the leading monomials of elements of $G$: Suppose
that $r\in I^*$ is such that $LM(r)|u$, and let $g\in G$ be such
that $LM(g)|LM(r)$: then $LM(g)|u$. 

\begin{proposition}[\cite{li}, Prop.~3.3, p.~70]
The set of monomials that are normal mod $I$ projects to 
a vector space basis for $T/I$. 
\end{proposition}

\end{document}